\documentclass{article}
\usepackage{amssymb}
\usepackage{graphicx}
\usepackage{amsmath}
\usepackage[singlespacing]{setspace}
\usepackage{fancyhdr}

\newtheorem{theorem}{Theorem}[section]

\newtheorem{lemma}[theorem]{Lemma}
\newtheorem{notation}[theorem]{Notation}

\newtheorem{proposition}[theorem]{Proposition}

\newenvironment{proof}[1][Proof]{\noindent\textbf{#1.} 
}{\ \rule{0.6em}{0.6em}}
\input tcilatex 
\begin{document}

\title{The Ratliff-Rush Closure of Initial Ideals of Certain Prime Ideals}
\author{Ibrahim Al-Ayyoub}
\maketitle

\begin{abstract}
Let $K$ be a field and let $m_{0},...,m_{n}$ be an almost arithmetic
sequence of positive integers. Let $C$ \ be a monomial curve in the affine $%
\left( n+1\right) $-space, defined parametically by $x_{0}=t^{m_{0}},\ldots
,x_{n}=t^{m_{n}}$. In this article we prove that the initial ideal of the
defining ideal of $C$\ is Ratliff-Rush closed.
\end{abstract}

\section*{Introduction\ \ \ \ \ \ \ \ \ \ \ \ \ \ \ \ \ \ \ \ }

In Section~\ref{RRclosure} we introduce the Ratliff-Rush closure of an ideal
and refer to some procedures used to compute it. In Section~\ref{MonoCurves}
we recall the Groebner bases of the prime ideals that are the defining
ideals of monomial curves as a result of a previous study. Section~\ref%
{MainResult} contains the main result of this article proving that the
initial ideals of these prime ideals are Ratliff-Rush closed.

\section{The Ratliff-Rush Closure\label{RRclosure}}

Let $R$ be a commutative Noetherian ring with unity and $I$ a regular ideal
in $R$, that is, an ideal that contains a nonzerodivisor. Then the ideals of
the form $I^{n+1}:I^{n}=\{x\in R\mid xI^{n}\subseteq I^{n+1}\}$ increase
with $n$. Let us denote%
\begin{equation*}
\widetilde{I}=\underset{n\geq 1}{\cup }(I^{n+1}:I^{n}).
\end{equation*}%
As $R$ is Noetherian, $\widetilde{I}$ $=I^{n+1}:I^{n}$ for all sufficiently
large $n$. Ratliff and Rush (1978) [Theorem 2.1] proved that $\widetilde{I}$%
\ is the unique largest ideal for which $(\widetilde{I})^{n}=I^{n}$ for
sufficiently large $n$. The ideal $\widetilde{I}$\ is called the Ratliff-Rush%
\textit{\ closure} of $I$ and $I$ is called \textit{Ratliff-Rush closed} if $%
I=\widetilde{I}$. It is easy to see that $I\subseteq \widetilde{I}$ and that
an element of $(I^{n}:I^{n+1})$ is integral over $I$. Hence for all regular
ideals $I$,%
\begin{equation*}
I\subseteq \widetilde{I}\subseteq \bar{I}\subseteq \sqrt{I}.
\end{equation*}%
where $\bar{I}$\ is the integral closure of $I$. Thus all radical and
integrally closed regular ideals are Ratliff-Rush closed. But there are many
ideals which are Ratliff-Rush closed but not integrally closed.

\ \ 

Rossi and Swanson (2003)\ examine the behavior of the Ratliff-Rush closure
with respect to some properties such as the Ratliff-Rush closure of powers
of ideals. They established new classes of ideals for which all the powers
are Ratliff-Rush closed. They also show that the Ratliff-Rush closure does
not behave well under several properties, such as, taking powers of ideals,
leading terms ideals, and the minimal number of generators. They present
many examples illustrating the different behaviors of the Ratliff-Rush
closure.

\ \ 

As yet, there is no algorithm to compute the Ratliff-Rush closure for
regular ideals in general. To compute $\cup _{n}(I^{n+1}:I^{n})$ we need to
find a positive integer $N$ such that $\cup _{n}(I^{n+1}:I^{n})$ $%
=I^{N+1}:I^{N}$. However, $I^{n+1}:I^{n}=I^{n+2}:I^{n+1}$\ does not imply
that $I^{n+1}:I^{n}=I^{n+3}:I^{n+2}$\ (see Example 1.8 in Rossi and Swanson
(2003)). Some different approaches have been used to decide the Ratliff-Rush
closure; Heinzer et al. (1992) established that a regular ideal $I$ (and
also every powers of I) is Ratliff-Rush closed if and only if the associated
graded ring, $gr_{I}(R)=\oplus _{n\geq 0}I^{n}/I^{n+1}$, has a
nonzerodivisor (has a positive depth). Elias (2003) established a procedure
for computing the Ratliff-Rush closure of $\mathbf{m}$-primary ideals of a
Cohen-Macaulay local ring with maximal ideal $\mathbf{m}$.

\ \ 

From the definition, it is clear that the Ratliff-Rush closure of a monomial
ideal is a monomial ideal, and this makes some computations easier. The
following two theorems and proposition serve us as a technique to compute
the Ratliff-Rush closure of the monomial ideals of interest in this article.

\begin{lemma}
Let $R,S$ be Noetherian rings. Assume $R$ is a faithfully flat $S$-algebra
and $I$ $\subset S$ an ideal. Then $IR$ is Ratliff-Rush closed in $R$ iff $I$
is Ratliff-Rush closed in $S$.
\end{lemma}

\begin{proposition}
\label{RRiffRR}Let $R=K[x_{0},...,x_{n}]$ and $S=K[x_{0},...,x_{m}]$ with $%
m\leq n$ where $K$ is a field. Let $I$ $\subset S$ be an ideal. Then $IR$ is
Ratliff-Rush closed in $R$ iff $I$ is Ratliff-Rush closed in $S$.
\end{proposition}

\begin{theorem}
\label{MainRR}Let $I$ be an ideal in the polynomial ring $%
R=K[x_{0},...,x_{n}]$ with $K$ a field. Let $r\geq 1$. If $I$\ is primary to 
$(x_{r},...,x_{n})$ and $\widetilde{I}\cap (I:(x_{r},...,x_{n}))\subseteq I$
\ then $I$\ is Ratliff-Rush closed.
\end{theorem}

\begin{proof}
Assume $I$ is not Ratliff-Rush closed. Let $m$ be an element such that $m\in 
\widetilde{I}$ $\backslash $ $I$. As $I$\ is primary to $(x_{r},...,x_{n})$
then there exists an integer $k$ such that $(x_{r},...,x_{n})^{k}\subseteq
I. $ In particular, $(x_{r},...,x_{n})^{l}m\subseteq I$ for some $l$. Choose 
$l\geq 1$ the smallest possible such integer. Then $(x_{r},...,x_{n})^{l-1}m%
\nsubseteqq I$. Let $m^{\prime }\in (x_{r},...,x_{n})^{l-1}$ be a monomial
such that $m^{\prime }m\notin I$. Then $(x_{r},...,x_{n})m^{\prime
}m\subseteq (x_{r},...,x_{n})^{l}m\subseteq I$. Thus $m^{\prime }m\in
I:(x_{r},...,x_{n})$ and $m^{\prime }m\in \widetilde{I}$ as $m\in \widetilde{%
I}$ . Therefore, $m^{\prime }m\in \widetilde{I}\cap
(I:(x_{r},...,x_{n}))\backslash I$.
\end{proof}

\section{{\protect\Large The Defining Ideals of Certain Monomial Curves 
\label{MonoCurves}}}

Let $n\geq 2$, $K$ a field and let $x_{0},...,x_{n},t$ be indeterminates.
Let $m_{0},...,m_{n}$ be an almost arithmetic sequence of positive integers,
that is, some $n-1$ of these form an arithmetic sequence, and assume $%
gcd(m_{0},...,m_{n})=1$. Let $P$\ be the kernel of the $K$-algebra
homomorphism $\eta :K[x_{0},...,x_{n}]\rightarrow K[t]$, defined by $\eta
(x_{i})=t^{m_{i}}$. A set of generators for the ideal $P$ was explicitly
constructed in Patil and Singh (1990). We call these generators the \textit{%
\textquotedblleft Patil-Singh generators\textquotedblright }. In a previous
study we proved that Patil-Singh generators form a Groebner basis for the
prime ideal $P$ with respect to the grevlex monomial order using the grading 
$wt(x_{i})=m_{i}$ with\textit{\ }$x_{0}<x_{1}<\cdots <x_{n}$ ( in this case 
\textit{\ }$\prod\limits_{i=0}^{n}x_{i}^{a_{i}}>_{grevlex}\prod%
\limits_{i=0}^{n}x_{i}^{b_{i}}$ if in the ordered tuple $%
(a_{1}-b_{1},...,a_{n}-b_{n})$\ the left-most nonzero entry is negative).
Before we state the Groebner basis we need to introduce some notations and
terminology that Patil and Singh (1990) used in their construction of the
generating set for the ideal $P$.

\ 

Let $n\geq 2$ be an integer and let $p=n-1$ . Let $m_{0},...,m_{p},m_{n}$ be
an almost arithmetic sequence of positive integers and $%
gcd(m_{0},...,m_{n})=1$, $0<m_{0}<\cdots <m_{p}$, and $m_{n}$ is arbitrary.
Let $\Gamma $ denote the numerical semigroup that is minimally generated by $%
m_{0},...,m_{p},m_{n}$, i.e. $\Gamma =\sum\limits_{i=0}^{n}\mathbb{N}_{%
\mathbf{0}}m_{i}$ . Put $\Gamma ^{\prime }=\sum\limits_{i=0}^{p}\mathbb{N}_{%
\mathbf{0}}m_{i}$\ and $\Gamma =\Gamma ^{\prime }+\mathbb{N}_{\mathbf{0}%
}m_{n}$.

\begin{notation}
\textit{For }$c,d\in \mathbb{Z}~$\textit{\ let }$[c,d]=\{t\in \mathbb{Z}\mid
c\leq t\leq d\}$\textit{. For }$t\geq 0$\textit{, let }$q_{t}\in \mathbb{Z}$%
\textit{, }$r_{t}\in \lbrack 1,p]$\textit{\ and }$g_{t}\in \Gamma ^{\prime }$%
\textit{\ \ be defined by }$t=q_{t}p+r_{t}$\textit{\ and }$%
g_{t}=q_{t}m_{p}+m_{r_{t}}$.
\end{notation}

Let $S=\{\gamma \in \Gamma \mid \gamma -m_{0}\notin \Gamma \}$. The
following is a part of Lemma (1.6) given in Patil (1993) that gives an
explicit description of $S$.

\begin{lemma}
(Patil (1993) Lemma 1.6))\textbf{\ }\textit{Let }$u=min\{t\geq 0\mid
g_{t}\notin S\}$\textit{\ and }$\upsilon =min\{b\geq 1\mid bm_{n}\in \Gamma
^{\prime }\}$\textit{. Then there exist unique integers }$w\in \lbrack
0,\upsilon -1]$\textit{, }$z\in \lbrack 0,u-1]$\textit{, }$\lambda \geq 1$%
\textit{, }$\mu \geq 0$\textit{, and }$\nu \geq 2$\textit{\ such that\newline
\ \ (i) }$g_{u}=\lambda m_{0}+wm_{n}$\textit{;\newline
\ \ (ii) }$\upsilon m_{n}=\mu m_{0}+g_{z}$\textit{;\newline
\ \ (iii) }$g_{u-z}+(\upsilon -w)m_{n}=\left\{ 
\begin{tabular}{ll}
$\left( \lambda +\mu +1\right) m_{0}\text{,}$ & if$\text{\ \ }r_{u-z}<r_{u}%
\text{;}$ \\ 
$\left( \lambda +\mu \right) m_{0}\text{,}$ & $\text{if \ }r_{u-z}\geq r_{u}%
\text{.}$%
\end{tabular}%
\right. $
\end{lemma}

\begin{notation}
Let\textit{\ }$q=q_{u},$\textit{\ }$r=r_{u}$\textit{. }For the rest of this
article the symbols $q,r,u,\upsilon ,w$, $z,\lambda $ and $\mu $ will have
the meaning assigned to them by the lemma and the notations above.
\end{notation}

Let $\ \varepsilon =\left\{ 
\begin{tabular}{ll}
$0\text{,}$ & if$\text{ \ }r>r_{z}\text{;}$ \\ 
$1\text{,}$ & if$\text{ \ }r\leq r_{z}$,%
\end{tabular}%
\right. $\newline
\ \ \newline
We state Patil-Singh generators as follows:\ \ \ \ \ \ \ \ \ \ \ \ 

\ \ \ \ \ \ \ \ \ \ \ \ \ \ \ \ \newline

\begin{tabular}{lll}
$\ \varphi _{i}$ & $=x_{i+r}x_{p}^{q}-x_{0}^{\lambda -1}x_{i}x_{n}^{w}$, & 
for $\ 0\leq i\leq p-r$; \\ 
$\psi _{j}$ & $=x_{\varepsilon p+r-r_{z}+j}x_{p}^{q-q_{z}-\varepsilon
}x_{n}^{\upsilon -w}-x_{0}^{\lambda +\mu -\varepsilon }x_{j}$, & for $\ j\in %
\left[ 0\ ,(1-\varepsilon )p+r_{z}-r\right] $; \\ 
$\theta $ & $=x_{n}^{\upsilon }-x_{0}^{\mu }x_{r_{z}}x_{p}^{q_{z}}$, &  \\ 
$\alpha _{i,j}$ & $=x_{i}x_{j}-x_{i-1}x_{j+1}$, & for $\ 1\leq i\leq j\leq
p-1$.%
\end{tabular}

\ \ \ \ \ 

\begin{theorem}
\label{GB}(Al-Ayyoub 2004))The set $\{\varphi _{i}\mid 0\leq i\leq p-r\}\cup
\{\theta \}$ $\cup $ $\{\alpha _{i,j}\mid 1\leq i\leq j\leq p-1\}$ $\cup $ $%
\{\psi _{j}\mid 0\ \leq \ j\leq (1-\varepsilon )p+r_{z}-r\}$ forms a
Groebner basis for the ideal $P$ with respect to the grevlex monomial order\
with $x_{0}<x_{1}<\cdots <x_{n}$ and with the grading $wt(x_{i})=m_{i}$.
\end{theorem}

\section{The Main Result\label{MainResult}}

In this section we prove that the initial ideal $inP$, of the defining ideal
of the monomial curves introduced in Section~\ref{MonoCurves}, is
Ratliff-Rush closed. The previous section states a Groebner basis for the
defining ideal $P$ with respect to the grevlex monomial order with the
grading $wt(x_{i})=m_{i}$ with $x_{0}<x_{1}<\cdots <x_{n}$. Therefore, $inP$
is generated by the following monomials\newline

\begin{tabular}{ll}
$x_{i}x_{p}^{q}$, & for \ $i$\ $\in \left[ r,p\right] $; \\ 
$x_{j}x_{p}^{q-q_{z}-\varepsilon }x_{n}^{\upsilon -w}$, & for $\ j\in \left[
\varepsilon p+r-r_{z},p\right] $; \\ 
$x_{n}^{\upsilon }$, &  \\ 
$x_{i}x_{j}$, & for\ \ \ $1\leq i\leq j\leq p-1$.%
\end{tabular}

\ \ \ 

Now we state the main result of the article:\ \ \ \ \ \ \ \ \ \ \ \ \ \ \ \ 

\begin{theorem}
\label{MainThm-RR}Let $P$ be the defining ideal of the monomial curves as
defined before. Then the ideal $inP$ is Ratliff-Rush closed.
\end{theorem}

Here is an outline for the proof of Theorem~\ref{MainThm-RR}: from the
generators above, it is clear that the monomial ideal $inP$ is primary to $%
(x_{1},...,x_{n})$. Therefore we can use Theorem~\ref{MainRR} to prove that $%
\left( inP\right) R$ is Ratliff-Rush closed in the polynomial ring $%
R=K[x_{1},...,x_{n}]$, and hence by Proposition~\ref{RRiffRR} Ratliff-Rush
closed in the polynomial ring $K[x_{0},...,x_{n}]$. In order to establish
the details of this outline we need to compute $(inP:(x_{1},...,x_{n}))/inP$%
. The following proposition is the first step in doing so.

\begin{proposition}
\label{x1---xp-1}with notation as before, then $%
(inP:(x_{1},...,x_{p-1}))/inP=(\overline{x}_{1},..$. ,$\overline{x}_{p-1})$,
where $\overline{x}_{i}$ is the image of $x_{i}$\ in the ring $R/inP$.
\end{proposition}

\begin{proof}
Let $\lambda =\min \{r,\varepsilon p+r-r_{z}\}$ and let $\sigma =\max
\{r,\varepsilon p+r-r_{z}\}$. Note that $(inP:(x_{i}))/inP=(\overline{x}%
_{1},...,\overline{x}_{p-1})$ for $1\leq i<\lambda $, and $%
(inP:(x_{i}))/inP=(\overline{x}_{1},...,\overline{x}_{p-1},\varepsilon
x_{p}^{q},(1-\varepsilon )\overline{x}_{p}^{q-q_{z}-\varepsilon }\overline{x}%
_{n}^{\upsilon -w})$ for $\lambda \leq i<\sigma $. Also note that $%
(inP:(x_{i}))/inP=(\overline{x}_{1},...,\overline{x}_{p-1},\overline{x}%
_{p}^{q},\overline{x}_{p}^{q-q_{z}-\varepsilon }\overline{x}_{n}^{\upsilon
-w})$ for $\sigma <i\leq p-1$.\thinspace Hence, it follows that\thinspace $%
(inP:(x_{1},...,x_{p-1}))/inP=\tbigcap\limits_{i=1}^{p-1}(inP:(x_{i}))/inP%
\,=\,(\overline{x}_{1},...,\overline{x}_{p-1})$.
\end{proof}

\begin{notation}
To simplify notations, in the sequel if a monomial happens to have an
indeterminate with a negative exponent then that monomial is treated as $0$.
For example, $x_{1}^{-2}x_{3}+x_{2}^{2}-x_{3}$ is $x_{2}^{2}-x_{3}.$\ 
\end{notation}

\begin{proposition}
\label{x1---xp}Let $p=n-1$ as before, then $(inP:(x_{1},...,x_{p}))/inP$ is
minimally generated in $K[x_{1},...,x_{n}]/inP$ by $\{\overline{x}_{i}%
\overline{x}_{p}^{q}\mid 1\leq i\leq r-1\}\cup \{\overline{x}_{i}\overline{x}%
_{p}^{q-1}\mid r\leq i\leq p-1\}\cup \{\overline{x}_{i}\overline{x}%
_{p}^{q-q_{z}-\varepsilon }\overline{x}_{n}^{\upsilon -w}\mid 1\leq i\leq
\varepsilon p+r-r_{z}-1\}\cup \{\overline{x}_{i}\overline{x}%
_{p}^{q-q_{z}-\varepsilon -1}\overline{x}_{n}^{\upsilon -w}\mid \varepsilon
p+r-r_{z}\leq i\leq p-1\}$.
\end{proposition}

\begin{proof}
We need to compute $\left(
\tbigcap\limits_{i=1}^{p-1}(inP:(x_{i}))/inP\right) \cap (inP:(x_{p}))/inP$.
Note that $(inP:(x_{p}))/inP$ is minimally generated by the following set of
monomials $\left\{ \overline{x}_{\varepsilon p+r-r_{z}}\overline{x}%
_{p}^{q-q_{z}-\varepsilon -1}\overline{x}_{n}^{\upsilon -w},\ldots ,%
\overline{x}_{p-1}\overline{x}_{p}^{q-q_{z}-\varepsilon -1}\overline{x}%
_{n}^{\upsilon -w},\overline{x}_{p}^{q-q_{z}-\varepsilon }\overline{x}%
_{n}^{\upsilon -w}\right\} \cup $\newline
$\left\{ \overline{x}_{r}\overline{x}_{p}^{q-1},\ldots \newline
,\overline{x}_{p-1}\overline{x}_{p}^{q-1},\overline{x}_{p}^{q}\right\} $. As
the intersection of two monomial ideals is generated by the least common
multiple of the monomial generators of each of the two ideals, then the
proposition follows by Proposition~\ref{x1---xp-1}.
\end{proof}

\ \ \ \ 

We next compute $(inP:(x_{n}))/inP$ . For the sake of notation we do so in
two cases. Also, at the same time we will prove Theorem~\ref{MainThm-RR} for
each of these cases separately. With the notations from Section~\ref%
{MonoCurves} consider the following two cases: Case 1: $\varepsilon >0$ or $%
q_{z}>0$, and Case 2: $\varepsilon =q_{z}=0$.

\subsection{Case 1: $\protect\varepsilon >0$ or $q_{z}>0$\ \ }

In this case $inP$ is generated by the following set of monomials\newline

\begin{tabular}{ll}
$x_{i}x_{p}^{q}$, & for \ $r\leq i\leq p$; \\ 
$x_{j}x_{p}^{q-q_{z}-\varepsilon }x_{n}^{\upsilon -w}$, & for $\ \varepsilon
p+r-r_{z}\leq $\ $j\leq p$; \\ 
$x_{n}^{\upsilon }$, &  \\ 
$x_{i}x_{j}$, & for\ $\ 1\leq i\leq j\leq p-1$.%
\end{tabular}
\newline
\ \ \ \ \newline
Therefore, $(inP:(x_{n}))/inP$ is minimally generated by \newline
$\{\overline{x}_{p+r-r_{z}}\overline{x}_{p}^{q-q_{z}-\varepsilon }\overline{x%
}_{n}^{\upsilon -w-1},\ldots ,\overline{x}_{p-1}\overline{x}%
_{p}^{q-q_{z}-\varepsilon }\overline{x}_{n}^{\upsilon -w-1}\}\cup \left\{ 
\overline{x}_{p}^{q-q_{z}-\varepsilon +1}\overline{x}_{n}^{\upsilon
-w-1}\right\} \cup \left\{ \overline{x}_{n}^{\upsilon -1}\right\} $. As the
intersection of two monomial ideals is generated by the least common
multiple of the monomial generators of each of the two ideals, then by
Proposition~\ref{x1---xp} it is straightforward to compute that $%
inP:((x_{1},...,x_{n}))/inP=(\tbigcap\limits_{i=1}^{n}(inP:(x_{i}))/inP=(%
\tbigcap\limits_{i=1}^{p}(inP:(x_{i}))/inP\cap (inP:(x_{n}))/inP$ is
generated by the monomials in the set $\varrho \cup \chi $, where $\varrho
=\{\overline{x}_{i}\overline{x}_{p}^{q-q_{z}-\varepsilon -1}\overline{x}%
_{n}^{\upsilon -1}\mid \varepsilon p+r-r_{z}\leq i\leq p-1\}$ and $\chi $
consists of the following monomials

\ \ 

\begin{tabular}{ll}
$\overline{x}_{i}\overline{x}_{p}^{q}\overline{x}_{n}^{\upsilon -w-1}$, & 
for $\ 1\leq i\leq r-1$; \\ 
$\delta _{q_{z}0}\overline{x}_{i}\overline{x}_{p}^{q-1}\overline{x}%
_{n}^{\upsilon -w-1}$, & for $\ r\leq i\leq \varepsilon p+r-r_{z}-1$; \\ 
$\overline{x}_{i}\overline{x}_{p}^{q-1}\overline{x}_{n}^{\upsilon -w-1}$, & 
for \ $\varepsilon p+r-\varepsilon r_{z}\leq i\leq p-1$;\  \\ 
$\overline{x}_{i}\overline{x}_{p}^{q-q_{z}-\varepsilon }\overline{x}%
_{n}^{\upsilon -1}$, & for \ $1\leq i\leq \varepsilon p+r-r_{z}-1$.%
\end{tabular}

\ \ \ 

Therefore, the preimages of the monomials in $\varrho \cup \chi $\ are the
only monomials in $(inP:(x_{1},...,x_{p}))\backslash inP$ $\ $in the ring$\
K[x_{1},...,x_{n}]$. By Theorem~\ref{MainRR} we prove that $inP$\ is
Ratliff-Rush closed by showing that none of these monomials belongs to the
Ratliff-Rush closure $\widetilde{inP}$ of $inP$. We show this separately for
the monomials in $\varrho $ and the monomials in $\chi $. First, assume $%
\overline{x}_{i}\overline{x}_{p}^{q-q_{z}-\varepsilon -1}\overline{x}%
_{n}^{\upsilon -1}\in \varrho $ is in $\widetilde{inP}$ for $\varepsilon
p+r-r_{z}\leq i\leq p-1$. Then by the definition of the Ratliff-Rush closure
we must have $\overline{x}_{i}\overline{x}_{p}^{q-q_{z}-\varepsilon -1}%
\overline{x}_{n}^{\upsilon -1}(x_{i}^{2})^{m}\in (inP)^{m+1}$ for some $%
m\geq 1$. By degree count for $x_{p}$ and $x_{n}$ we must have $\overline{x}%
_{i}\overline{x}_{p}^{q-q_{z}-\varepsilon -1}\overline{x}_{n}^{\upsilon
-1}(x_{i}^{2})^{m}\in (x_{i}^{2})^{m+1}$, contradiction by the $x_{i}$
degree count.

\ \ \ 

Now assume $x_{i}x_{p}^{a}x_{n}^{b}$ is a monomial in $\chi $ ($a\leq q$ and 
$b<\upsilon $) such that $x_{i}x_{p}^{a}x_{n}^{b}$ $\in \widetilde{inP}$.
Then $x_{i}x_{p}^{a}x_{n}^{b}$ $(x_{i}^{2})^{m}\in (inP)^{m+1}$ for some $%
m\geq 1$. By $x_{n}$ and $x_{i}$-degree count for $1\leq i\leq p-1$ we must
have $x_{i}^{2m+1}x_{p}^{a}x_{n}^{b}\in $($\delta _{i\geq
r}x_{i}x_{p}^{q},\delta _{i\geq \varepsilon
p+r-r_{z}}x_{i}x_{p}^{q-q_{z}-\varepsilon }$ $x_{n}^{\upsilon -w}$)$^{m+1}$.
Note if $a=q$ then we must have $i<r$, thus $x_{i}^{2m+1}x_{p}^{a}x_{n}^{b}%
\in (\delta _{i\geq \varepsilon p+r-r_{z}}x_{i}$ $x_{p}^{q-q_{z}-\varepsilon
}x_{n}^{\upsilon -w})^{m+1}$. Assume $a<q$. Then $%
x_{i}^{2m+1}x_{p}^{a}x_{n}^{b}\notin (\delta _{i\geq r}x_{i}x_{p}^{q})$,
hence $x_{i}^{2m+1}x_{p}^{a}x_{n}^{b}$ $\in $ ($\delta _{i\geq \varepsilon
p+r-r_{z}}x_{i}x_{p}^{q-q_{z}-\varepsilon }x_{n}^{\upsilon -w})^{m+1}$. In
either case it implies that implies $i\geq \varepsilon p+r-r_{z}$ and $b\geq
\upsilon -w$. But there are no such monomials in $\chi $.

\subsection{Case 2: $\protect\varepsilon =q_{z}=0$}

In this case $inP$ is minimally generated by the following set of monomials%
\newline

\begin{tabular}{ll}
$x_{i}x_{p}^{q}$, & for \ $r\leq i\leq p$; \\ 
$x_{j}x_{p}^{q}x_{n}^{\upsilon -w}$, & for $\ r-r_{z}\leq $\ $j\leq r-1$; \\ 
$x_{n}^{\upsilon }$, &  \\ 
$x_{i}x_{j}$, & for$\ \ 1\leq i\leq j\leq p-1$.%
\end{tabular}

\ \ \ \ \ \ \ \ \ \ \ \newline
Therefore, $(inP:(x_{n}))/inP$ is minimally generated by \newline
$\{\overline{x}_{r-r_{z}}\overline{x}_{p}^{q}\overline{x}_{n}^{\upsilon
-w-1},\ldots ,\overline{x}_{r-1}\overline{x}_{p}^{q}\overline{x}%
_{n}^{\upsilon -w-1}\}\cup \left\{ \overline{x}_{n}^{\upsilon -1}\right\} $.
By Proposition~\ref{x1---xp} it follows that $inP:((x_{1},...,x_{n}))/inP=(%
\tbigcap\limits_{i=1}^{n}(inP:(x_{i}))/inP=(\tbigcap%
\limits_{i=1}^{p}(inP:(x_{i}))/inP\cap (inP:(x_{n}))/inP$ is generated by
the monomials in the set $\varrho \cup \chi $, where $\varrho =\{\overline{x}%
_{i}\overline{x}_{p}^{q-1}\overline{x}_{n}^{\upsilon -1}\mid r-r_{z}\leq
i\leq p-1\}$ and $\chi $ consists of the following monomials

\ 

\begin{tabular}{ll}
$\overline{x}_{i}\overline{x}_{p}^{q}\overline{x}_{n}^{\upsilon -1}$, & for
\ $1\leq i\leq r-r_{z}-1$; \\ 
$\overline{x}_{i}\overline{x}_{p}^{q}\overline{x}_{n}^{\upsilon -w-1}$, & 
for $\ r-r_{z}\leq i\leq r-1$;%
\end{tabular}%
\newline
\ \ \newline
Therefore, the preimages of the monomials in $\varrho \cup \chi $\ are the
only monomials in $(inP:(x_{1},...,x_{p}))\backslash inP\ $in the ring$\
K[x_{1},...,x_{n}]$. By Theorem~\ref{MainRR} we prove that $inP$\ is
Ratliff-Rush closed by showing that none of these monomials belongs to the
Ratliff-Rush closure $\widetilde{inP}$ of $inP$. We show this separately for
the monomials in $\varrho $ and the monomials in $\chi $. First, assume $%
\overline{x}_{i}\overline{x}_{p}^{q-1}\overline{x}_{n}^{\upsilon -1}\in
\varrho $ is in $\widetilde{inP}$ for $r-r_{z}\leq i\leq p-1$. Then by the
definition of the Ratliff-Rush closure we must have $\overline{x}_{i}%
\overline{x}_{p}^{q-q_{z}-\varepsilon -1}\overline{x}_{n}^{\upsilon
-1}(x_{i}^{2})^{m}\in (inP)^{m+1}$ for some $m\geq 1$. By degree count for $%
x_{p}$ and $x_{n}$ we must have $\overline{x}_{i}\overline{x}_{p}^{q-1}%
\overline{x}_{n}^{\upsilon -1}(x_{i}^{2})^{m}\in (x_{i}^{2})^{m+1}$,
contradiction by the $x_{i}$ degree count.

\ \ \ 

Now assume $x_{i}x_{p}^{q}x_{n}^{b}$ is a monomial in $\chi $ ( $b<\upsilon $%
) such that $x_{i}x_{p}^{q}x_{n}^{b}$ $\in \widetilde{inP}$. Then $%
x_{i}x_{p}^{q}x_{n}^{b}$ $(x_{i}^{2})^{m}\in (inP)^{m+1}$ for some $m\geq 1$%
. By $x_{n}$ and $x_{i}$-degree count for $1\leq i\leq p-1$ we must have $%
x_{i}^{2m+1}x_{p}^{q}x_{n}^{b}\in (\delta _{i\geq r}x_{i}x_{p}^{q},\delta
_{i\geq r-r_{z}}x_{i}$ $x_{p}^{q}x_{n}^{\upsilon -w})^{m+1}$. Note we must
have $i<r$, thus $x_{i}^{2m+1}x_{p}^{q}x_{n}^{b}\in (\delta _{r-r_{z}\leq
i\leq r-1}x_{i}$ $x_{p}^{q}x_{n}^{\upsilon -w})^{m+1}$. This implies $%
r-r_{z}\leq i\leq r-1$ and $b\geq \upsilon -w$. But there are no such
monomials in $\chi $.

\section*{Acknowledgement}

The author thanks Prof. Swanson I. for the useful discussions and comments
during the course of this work.

{\small Department of Mathematics and Statistics}

{\small Jordan University of Science and Technology}

{\small P O Box 3030, Irbid 22110, Jordan}

{\small Email address: iayyoub@just.edu.jo}

\end{document}